\documentclass{article}
\usepackage[numbers]{natbib}
\newcommand{\mainmatter}{}
\newcommand{\titlerunning}[1]{}
\newcommand{\authorrunning}[1]{}
\newcommand{\email}[1]{\texttt{#1}}

\makeatletter
\def\@maketitle{%
  \newpage
  \null
  \vskip 2em
  \begin{center}%
  \let \footer \@empty
    {\large \bfseries \@title \par}
    \vskip 1.5em
    {\small                        
      \lineskip .5em%
      \begin{tabular}[t]{c}%
        \@author
      \end{tabular}\par}%
    \vskip 1em
    {\small \@date}
  \end{center}%
  \par
  \vskip 1.5em}                    
\makeatother

\newcommand{\institute}[1]{\date{#1}}

	\providecommand{\pgfsyspdfmark}[3]{} 

\label{packages}
	\usepackage{microtype}
	\usepackage[utf8]{inputenc}
	\usepackage{amsmath}
	\usepackage{amssymb}
	\usepackage{algorithm}
	\usepackage{mathtools}
	\usepackage{aligned-overset}
	\usepackage{blindtext}
	\usepackage{booktabs} 
	\usepackage{csquotes}
	\usepackage[dvipsnames]{xcolor}
	\usepackage{enumitem}
		\setlist[enumerate]{itemsep=0.1em,topsep=0.4em,leftmargin=*,labelindent=1em}		
		\setlist[itemize]{itemsep=0.1em,topsep=0.4em,leftmargin=*,labelindent=1em,label={$\bullet$}}		
	\usepackage{environ}
	\usepackage{matlab-prettifier}
	\usepackage{units}
	\usepackage{resizegather}
	\usepackage{subcaption} 
		
	\usepackage{standalone} 
	\usepackage{tabularx}
	\usepackage{arydshln} 
	\usepackage{lipsum}

	\usepackage{matlab-prettifier}
	\usepackage{graphicx}  
		\definecolor{udsblue}    {RGB}{0,72,119}
		\definecolor{udsdarkblue}{RGB}{1,40,63}
		\definecolor{slightgray}{RGB}{215,215,215}
		\definecolor{whitegray}{rgb}{0.985,0.985,0.985}

\label{other}
	\newif\ifAppendix
	\newif\ifTodos   
	\newif\ifFormReview        
	\newif\ifIFACClassFile
	\newif\ifOldIEEEClassFile 
	\newif\ifAcronymsVisible
	\newif\ifAnimation
	\newif\ifTwoColumn
	\newif\ifAcronymCount

\label{standard macros/user-defined commands}
	\newcounter{RemarkCounter}
	\usepackage{listings}
		\newcommand{\inv}{^{\raisebox{-0.0ex}{$\scriptscriptstyle - 1$}}}
		\newcommand{\tr}{^{\raisebox{0.37ex}{$\scriptscriptstyle\intercal$}}}

	\newcommand{\diag}{\operatorname{diag}}
	\DeclareMathOperator*{\argmin}{arg\,min}
	\renewcommand{\(}{\left(}
	\renewcommand{\)}{\right)}
	\newcommand{\numberthis}{\refstepcounter{equation}\tag{\theequation}}

	\Appendixfalse
	\TwoColumnfalse
	\IFACClassFilefalse
	\input{./BuildFiles/header_hyperref.tex}
	
	\usepackage{autonum}
\usepackage[nameinlink,noabbrev,capitalise]{cleveref}
	\def\case#1{\uppercase{#1}}
	
	\crefname{subsection}{\case{s}ubsection}{\case{s}ubsections}
	\crefname{table}{Table}{Tables}
	
	\usepackage{caption}
	\usepackage{float}   
	\newfloat{FloatList}{htbp}{lol}
	\floatname{FloatList}{List}           
	\crefname{FloatList}{List}{List}
	
	\newcounter{nonfloatlist}
	
	\crefname{nonfloatlist}{List}{List}

		\usepackage{mdframed}
	\usepackage[breakable]{tcolorbox}

	\newenvironment{bThm}[2][]{%
		\def\envPH{#2}%
		\begin{\envPH}[#1]
	}{%
		\end{\envPH}
	}

\usepackage{amsthm}

		{%
			\begin{algorithm}[#2]
			\caption{#1}
		}{%
			\end{algorithm}
		}

	\AcronymCounttrue
	\label{prepare acronyms (acc, accp)}
\usepackage[nolist,nohyperlinks,printonlyused]{acronym}

\definecolor{colAcro}{RGB}{0, 110, 14} 	
\ifAcronymsVisible
		\newcommand{\highlight}[1]{\textcolor{colAcro}{#1}}

		\let\oldacl\acl	  \renewcommand*\acl[1]{\highlight{\oldacl{#1}}}
		\let\oldacp\acp	  \renewcommand*\acp[1]{\highlight{\oldacp{#1}}}
		\let\oldaclp\aclp \renewcommand*\aclp[1]{\highlight{\oldaclp{#1}}}
		\let\oldAcl\Acl	  \renewcommand*\Acl[1]{\highlight{\oldAcl{#1}}}
		\let\oldAcp\Acp	  \renewcommand*\Acp[1]{\highlight{\oldAcp{#1}}}
		\let\oldAclp\Aclp \renewcommand*\Aclp[1]{\highlight{\oldAclp{#1}}}
\else
	\newcommand{\highlight}[1]{#1}
\fi

\newcommand{\contentAcronymUsageTable}{}

	\newcommand{\incrementacronym}[1]{%
	    \stepcounter{acronymcount#1}%
	}

		\makeatletter
		\newcommand{\noLongAcronym}[2]{%
		    \if\relax\detokenize{#1}\relax
		        #2
		    \else
		        \@noLongAcronym#1\@nil
		    \fi
		}
		\def\@noLongAcronym#1#2\@nil{%
		    \if l#1%
		    \else
		        #2
		    \fi
		}
		\makeatother

	\makeatletter
	\newcommand{\myacro}[3][]{%
		\edef\@ctrname{acronymcount#2}%
		\expandafter\newcounter\expandafter{\@ctrname}%
		\expandafter\gdef\csname increment#2\endcsname{%
			\incrementacronym{#2}%
		}%
		\acro{#2}{#3}%
		\gappto{\contentAcronymUsageTable}{
		 	\displayacronymcount[#1]{#2}%
		}%
	}
	
	\newcommand{\myacroWO}[3][]{
		\edef\@ctrname{acronymcount#2}%
		\expandafter\newcounter\expandafter{\@ctrname}%
		\expandafter\setcounter\expandafter{\@ctrname}{-1}
		\expandafter\gdef\csname increment#2\endcsname{%
		}%
		\acro{#2}{#3}%
		\gappto{\contentAcronymUsageTable}{
		 	\displayacronymcount[#1]{#2}%
		}%
	}
	
	\newcommand{\condWO}[4]{%
		\ifnum\csname c@acronymcount#2\endcsname=-1\relax
		   	\csname#3l#4\endcsname{#2}%
		\else
			\csname#3#1#4\endcsname{#2}%
		\fi%
		\noLongAcronym{#1}{\csname increment#2\endcsname}
	}
	
	\newcommand{\acc}[2][]{%
		\condWO{#1}{#2}{ac}{}
	}
	\newcommand{\accp}[2][]{%
		\condWO{#1}{#2}{ac}{p}
	}

	\makeatother

	\newcommand{\displayacronymcount}[2][]{%
		#2
		&\edef\temp{\noexpand\arabic{acronymcount#2}}\temp
	   	& #1
   		\\
	}

	\Todostrue
	\definecolor{burntorange}{rgb}{0.8, 0.33, 0.0}
\ifFormReview
	\newcommand{\form}[1]{\textcolor{burntorange}{#1}}

\else
	\newcommand{\form}[1]{#1}

\fi

\let\oldLipsum\lipsum
	\renewcommand{\lipsum}[1][]{\textcolor{lightgray}{\oldLipsum[#1]}}

	\usepackage[titles]{tocloft}
	\makeatletter
		\renewcommand{\@pnumwidth}{4em}
		\renewcommand{\@tocrmarg}{4em}
	\makeatother
	\newcommand{\Rset}{\mathbb{R}}
\newcommand{\uc}[1]{\expandafter\MakeUppercase\expandafter{#1}}

\newcommand{\ii}{i}
\newcommand{\jj}{j}
\newcommand{\kk}{k}
\renewcommand{\ll}{n}
\newcommand{\mm}{m}
\newcommand{\nCOS}{n}%

\newcommand{\SO}{\mathrm{SO}}

\newcommand{\NP}{\mathrm{NP}}
\newcommand{\UP}{\mathrm{UP}}

\newcommand{\cNPUP}{C_{\NP,\UP}}

\newcommand{\J}{\mathcal{J}}
\newcommand{\JMO}{\J}

\newcommand{\xMO}{x}

\newcommand{\xComp}{\xMO^{\prime}}

\newcommand{\param}{\hat}
\newcommand{\realize}{\check}
\newcommand{\norm}{\bar}
\newcommand{\paramAccent}{hat}
\newcommand{\realizeAccent}{check}

\newcommand{\stan}{}
\renewcommand{\ne}{\breve}
\newcommand{\neAcc}{breve}

\newcommand{\spanVecMat}{V}

		\newcommand{\costFun}{J}
		\newcommand{\costFunVec}{\mathcal{J}}
		\newcommand{\nObj}{{n_J}}
		\newcommand{\decSetFeas}{\mathbf{X}}
		\newcommand{\decSetPareto}{\mathbf{X}_{\mathrm{P}}}
		
		\newcommand{\imSetPareto}[1][]{#1{\mathbf{J}}_{\mathrm{P}}}

\newenvironment{mat}{
  \begin{bmatrix}
}{
  \end{bmatrix}
}

\NewDocumentCommand{\matInline}{O{} O{} m}{%
	#1[#3#2]
}

	\begin{acronym}
	
	\myacroWO{MO}{multi-objective}
	\myacroWO{OP}{optimization problem}
	\myacro{MOOP}{multi-objective optimization problem}
	\myacroWO{MOO}{multi-objective optimization}
	\myacro{pOP}{parameterized optimization problem}	
	\myacro{MPC}{model predictive control}
	\myacroWO{MOMPC}{multi-objective model predictive control}
	
	\myacro{PF}{Pareto front}
	\myacroWO{POS}{Pareto optimal solution}
	
	\myacro{NP}{nadir point}
	\myacro{UP}{utopia point}
	\myacro{SO}{shooting origin}
	\myacroWO{SDV}{shooting direction vector}
	
	\myacro{CHIM}{convex hull of individual minima}
	\myacro{IM}{individual minimum}
		\acroplural{IM}[IMs]{individual minima}
	\myacroWO{KP}{knee-point}
	
	\myacro{WS}{weighted-sum}
	\myacro{NBI}{normal boundary intersection}
	\myacro{SRI}{spherical ray intersection}
	\myacroWO{PS}{Pascoletti-Serafini}

	\myacroWO{DM}{decision-making}
	
	\myacroWO{HVAC}{heating, ventilation and air conditioning}
	\myacro{MRS}{marginal rate of substitution}

\end{acronym}

	\mainmatter
\title{%
	Non-Extreme Individual Minima
	for Improved Pareto Front Sampling Efficiency 
	and 
	Decision-Making
}
\titlerunning{Non-Extreme Individual Minima} 

\author{
	Markus Herrmann-Wicklmayr
	\and
	Kathrin Flaßkamp
} 
\authorrunning{Herrmann-Wicklmayr and Flaßkamp} 
\institute{%
    Systems Modeling and Simulation, Saarland University, Saarbrücken, Germany\\
    \email{
	    \{markus.herrmannwicklmayr,kathrin.flasskamp\}%
	    @uni-saarland.de
    }
}

	\graphicspath{{./figures/}}

	\theoremstyle{definition} 
	
	\crefname{definition}{Definition}{Definitions}

\newcommand{\constPE}{M}%
\newcommand{\constPPE}{L}%

\begin{document}

\maketitle

\begin{abstract}
	In multi-objective optimization, the set of optimal trade-\linebreak offs---the Pareto front---often contains regions that are extremely steep or flat.
	The Pareto optimal points in these regions are typically of limited interest for decision-making, as the marginal rate of substitution is extreme: 
	a marginal improvement in one objective necessitates a significant deterioration in at least one other objective.
	These unfavorable trade-offs frequently occur near the individual minima, where single objectives attain their minimum values without considering the remaining criteria.

	To address this, we propose the concept of \emph{non-extreme individual minima} that relies on the notion of $\constPPE$-practical proper efficiency.
	These points can serve as a less sensitive replacement for \emph{standard} individual minima in subsequent related methods.
	Specifically, they allow for a more practical restriction of the Pareto front sampling within a refined utopia-nadir hyperbox, 
	provide a meaningful basis for image space normalization, 
	and can enhance decision-making techniques, such as knee-point methods, by focusing on regions with acceptable trade-offs.

	We provide a computationally efficient algorithm to determine these non-extreme individual minima by solving at most $2n_J$ standard weighted-sum scalarizations, where $n_J$ is the number of objectives. 
	To ensure robustness across varying objective scales, the method incorporates an integrated image space normalization strategy. 
	Numerical examples, specifically a convex academic case and a non-convex real-world application,
	demonstrate that the method successfully excludes \form{practically irrelevant} regions in the image space.
	\\[0.9\baselineskip]
		
\end{abstract}


\section{Introduction}
In \accp{MOOP}, a solution 
is typically selected among all optimal trade-offs \cite{roy_multicriteria_1996}, i.e.\ the Pareto optimal points forming the so-called \acc{PF}.
This task must be performed either by a human decision-maker or by an automated \acc{DM} scheme.
In both cases, it is reasonable to exclude parts of the \acc{PF} that are extremely steep or flat, i.e. parts of the \acc{PF} are avoided where a minor improvement in one objective yields a major deterioration in at least one other component.
We denote the regions of the \acc{PF} that are extremely steep or flat as \form{\emph{practically irrelevant}}.

Conventionally, the identification of \form{\emph{relevant}} solutions is performed a-posteriori: 
a human decision-maker localizes preferred regions only after a fine and ideally uniform sampling of the \acc{PF} has been generated, which is a challenging task in itself \cite{eichfelder_adaptive_2008}. 
However, this manual procedure is severely limited by human perception; 
once the number of objectives exceeds three, the resulting data points become nearly impossible to visualize and interpret.

Furthermore, modern applications with real-time requirements often preclude human intervention entirely, making this a-posteriori approach infeasible even for low-dimensional problems. 
Such scenarios necessitate automated \acc{DM} methods. 
If these methods are guided by information about the \accp{IM} \cite{zavala_stability_2012,chiu_minimum_2016,herrmann-wicklmayr_individual_2025}, 
it becomes particularly reasonable to utilize more robust reference points.

To address these challenges, we propose an a-priori strategy based on the notion of $\constPPE$-practical proper efficiency. 
By utilizing what we denote as \emph{non-extreme} \accp{IM}, we can delimit the search space before the full sampling process begins. 
The core idea is to specifically target and remove unpromising regions in the proximity of the \enquote{standard} \accp{IM} where the \acc{MRS} is excessively high. 
Since these non-extreme \accp{IM} can be determined efficiently by solving at most $2n_J$ \acc{WS} scalarizations, they enable a more focused allocation of computational resources by avoiding \form{practically irrelevant} trade-offs from the outset. 
Beyond computational efficiency, this approach provides automated \acc{DM} schemes 
(as e.g. \cite{li_knee_2020,schmitt_incorporating_2022}) 
with less sensitive and more representative reference points, which ultimately increases the reliability of the resulting selection.

The paper is organized as follows:
In \cref{sec:preliminaries} we introduce the basic notation and definitions of \accp{MOOP}.
Furthermore, we define the term non-extreme \acc{IM}.
In \cref{sec:Derivation_Method} we show an intuitive way to compute non-extreme \accp{IM}. 
The method is then summarized in an algorithm.
In \cref{sec:Numerical_Example} we apply the method to an exemplary \acc{MOOP}.
Finally, the paper is concluded in \cref{sec:Conclusion_and_Outlook}.

\textbf{Notation.}
We denote the standard basis of $\Rset^{\nCOS}$ by $\{e_{1},\dots,e_{\nCOS}\}$.
For a vector $v \in \Rset^n$, the relations $v = 0$ and $v \geq 0$ are to be understood component-wise, 
i.e., $v_{i} = 0$ and $v_{i} \geq 0$ for all $i = 1, \dots, n$.
We use the operator $\diag(v)$, which returns a diagonal matrix with $v$ on its diagonal
and the sign function $\operatorname{sign}(x)$, defined as $\operatorname{sign}(x) = 1$ if $x > 0$ and $\operatorname{sign}(x) = -1$ otherwise.
Bold symbols $\boldsymbol{0}$ and $\boldsymbol{1}$ denote vectors of zeros and ones, respectively, with context-dependent dimensions.

\section{Preliminaries}
\label{sec:preliminaries}
We repeat the problem setting and the characterization of important quantities in \accp{MOOP} from \cite{herrmann-wicklmayr_individual_2025} in the next two subsections.

\subsection{Problem Statement}
Consider the \acc{MOOP} 
\begin{equation}
	\min_{\xMO \in \decSetFeas} \
	\costFunVec(\xMO)
	\tag{$\mathcal{P}$}
	\label{eq:MOOP}
\end{equation}
with the vector 
$\costFunVec(\xMO) := \matInline[\left][\right]{\costFun_1(\xMO), \ldots, \costFun_\nObj(\xMO)}\tr$
of $n_J$ objectives $\costFun_{\ii}:\decSetFeas \rightarrow \mathbb{R}$, $\ii=1,\ldots,n_J$.
The set $\decSetFeas \subseteq \mathbb{R}^{n_{\xMO}}$ denotes the feasible set; a point $x\in \mathbb{R}^{n_{\xMO}}$ is feasible if $x\in \decSetFeas$.
The vector-valued minimization in \eqref{eq:MOOP} is clarified by definitions and conventions adopted from \cite{stieler_performance_2018}.
\begin{bThm}[Pareto optimality, nondominance]{definition}
	\label{def:Nondom}
	A point $\xMO^{\star} \in \decSetFeas$ is an \textup{efficient} or a \textup{Pareto optimal} solution
	to the \acc{MO} \acc[l]{OP} \eqref{eq:MOOP} if there does not exist any feasible $\xMO \in \decSetFeas$ such that
	\begin{align}
		& \costFun_{\ii}(\xMO) \leq \costFun_{\ii}\left(\xMO^{\star}\right)
		\text{ for all }  \ii \in\{1, \ldots, \nObj\} 
		\text{ and }
		\\
		& \costFun_k(\xMO)<\costFun_k\left(\xMO^{\star}\right) 
		\text { for at least one } k \in\{1, \ldots, \nObj\}.
	\end{align}
	The respective image value 
	$
	\costFunVec\left(\xMO^{\star}\right)
	$ 
	is called \textup{nondominated}.
	The set of all nondominated points is the nondominated set or \textup{\acc[l]{PF}} $\imSetPareto:=\left\{\costFunVec(\xMO)  \ \middle| \ \xMO \in \decSetPareto\right\}$,
	with the Pareto set $\decSetPareto$ given by
	\begin{equation}	
		\decSetPareto := \argmin_{\xMO \in \decSetFeas} \, \costFunVec(\xMO)	
		= \left\{
			\xMO \in \decSetFeas 
		\ \middle| \ 
			\xMO \text{ is a \acc[s]{POS} to \eqref{eq:MOOP}} 
		\right\}.
	\end{equation}
\end{bThm}

\subsection{Characteristic Quantities}
We derive characteristic quantities of an \acc{MOOP}.
We denote by $\xMO_{\ii}^*$, $\ii = 1,\ldots,n_J$ the \accp[f]{IM}, i.e.~solutions to the single-objective \acc{OP} 
$
	\min_{\xMO \in \decSetFeas} \ J_{\ii}(\xMO), \ \ii=1,\dots,n_J.
$
Evaluating the full objective vector $\JMO$ at all \acc{IM} defines the pay-off matrix (see, e.g.~\cite{ehrgott_multicriteria_2005})
$\Phi=\matInline{\JMO(\xMO_{1}^*), \ldots, \JMO(\xMO_{n_J}^*)}$.

\begin{bThm}[Utopia and nadir point]{definition}
	The \textup{\acc[f]{UP}} $\J_\UP$ is defined as the row-wise minimum of $\Phi$.
	We define the \textup{\acc[f]{NP}} $\J_\NP$ as the row-wise maximum of $\Phi$.
	\label{defi:NP_UP}
\end{bThm}
Note that, in literature, it is sometimes further distinguished between pseudo \accp{NP} (as defined in \cref{defi:NP_UP}) and real \accp{NP} \cite[Chapter 2.2]{ehrgott_multicriteria_2005}.

\begin{bThm}[Normalized image space]{definition}
	Let $\cNPUP$ denote the positive definite, diagonal matrix
	$
		\cNPUP := 
		\diag(\J_{\NP}-\J_{\UP})\inv.
	$
	Then, the operation 
	$
		\norm{\J}(x) = \cNPUP (\J - \J_{\UP})
	$
	shifts and scales the image space, such that the \acc{PF} is contained in the unit box (or hyperbox) spanned by the two points $\boldsymbol{0}$ and $\boldsymbol{1}$. 
	We call this a \textup{normalized image space}.
	The bar accent indicates a quantity in that normalized image space, 
	e.g. $\bar{\J}_{\UP}=\boldsymbol{0}$ and $\bar{\J}_{\NP}=\boldsymbol{1}$ are the \acc{UP} and \acc{NP} in that normalized space.
\end{bThm}

\subsection{Non-Extreme Individual Minima}
\begin{figure}[h]
	\centering
	\includegraphics[width=0.45\linewidth,page=2]{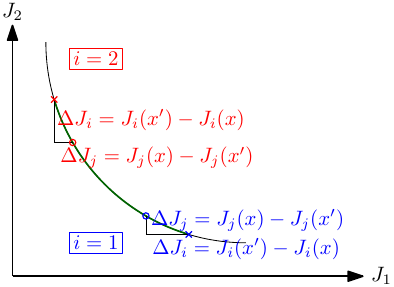}
	\caption{
		Objective value evaluated at $\xComp$ marked with a cross, which was obtained using the weight vector $w=(w_1,w_2)\tr$, and at $\xMO$ marked with a circle.
	}
	\label{fig:L_practical_PE}
\end{figure}

In order to define what we mean by \emph{non-extreme \accp{IM}} we first look at the definition of proper efficiency:
\begin{bThm}[Proper efficiency (in the sense of Geoffrion \cite{geoffrion_proper_1968})]{definition}
	\label{defi:proper_efficiency}
	A feasible point $\xComp \in \decSetFeas$ is a \textup{properly efficient} solution if it is efficient and if there exists some real number $\constPE>0$ such that, for each $\ii \in \{1,\ldots,n_J\}$ and $\xMO \in \decSetFeas$ satisfying $J_{\ii}(\xMO)<J_{\ii}(\xComp)$, there exists at least one $\jj \in \{1,\ldots,n_J\}\setminus \ii$ such that $J_{\jj}(\xComp)<J_{\jj}(\xMO)$ and	
	\begin{equation}
		\frac{J_{\ii}(\xComp)-J_{\ii}(\xMO)}{J_{\jj}(\xMO)-J_{\jj}(\xComp)} \leq \constPE.
	\end{equation}
\end{bThm}
We can illustrate \cref{defi:proper_efficiency} with an exemplary bi-objective \acc{OP} in \cref{fig:L_practical_PE}
where we investigate the proper efficiency of the solutions $\xComp$ marked with a cross.
For this it is helpful to consider a comparative solution $\xMO$, marked with a circle.
In the following we assume that $w_{\ii}>0, \ii \in \{1,\ldots,n_J\}$.
\\
We first fix the case where $\ii=2$ (and subsequently $\jj$ must be $1$).
For a convex \acc{PF} as in \cref{fig:L_practical_PE} we know that at the solution $\J(\xComp)$, 
that was obtained by solving the \acc{WS} problem $\min_{\xMO \in \decSetFeas} w\tr \J$, 
the \acc{PF} has the slope $m=-w_1/w_2$.
We notice that with $\xMO \rightarrow \xComp$ 
the approximation
\begin{equation}
	-w_1/w_2 
	= m 
	\approx \frac{
		-\(J_2(\xComp)-J_2(\xMO)\)
	}{
		J_1(\xMO)-J_1(\xComp)
	} 
	=: \frac{-\Delta J_2}{\Delta J_1}
\end{equation}
becomes exact.
We then obtain $w_1/w_2 = \lim_{\xMO \rightarrow \xComp} 
\Delta J_2 / \Delta J_1$.
In the same way, for the case $\ii=1$ (and $\jj=2$), we can derive $w_2/w_1 = \lim_{\xMO \rightarrow \xComp}
\Delta J_1 / \Delta J_2$.
This limit has a significant interpretation in \acc{MOO}: 
it represents the \acc[f]{MRS}. 
The \acc{MRS} quantifies the rate at which a decision-maker is willing to sacrifice an amount of objective $J_{\jj}$ to obtain a marginal improvement in objective $J_{\ii}$ while remaining on the \acc{PF}. 
In this context, the constant $\constPE$ in \cref{defi:proper_efficiency} acts as an upper bound on the \acc{MRS}.
This means that for both cases the solutions $\xComp$ are 
properly efficient with $\constPE=w_1/w_2$ and $\constPE=w_2/w_1$, respectively.
For the general $n_J$-dimensional case, it was proved in \cite{geromel_upper_1991,karimi_linear_2017} that 
\begin{equation}
	\frac{
		J_{\ii}(\xComp) - J_{\ii}(\xMO)
	}{
		J_{\jj}(\xMO) - J_{\jj}(\xComp)
	}
	\leq \frac{w_{\jj}}{w_{\ii}}, \quad \ii,\jj\in \{1,\ldots,n_J\}, \ \ii\neq \jj.
\end{equation}
Although \cref{defi:proper_efficiency} requires $\constPE$ to be finite, the value for it can still be arbitrarily high, implying an extreme trade-off where a negligible gain in one objective requires a massive loss in another. 
\form{%
From an engineering or economic perspective, solutions with an excessively high \acc{MRS} are not of practical relevance.
}%
In order to upper bound the value of $\constPE$ and thus the allowable \acc{MRS}, we introduce the following definition:
\begin{bThm}[$\constPPE$-practical proper efficiency]{definition}
	A feasible point $\xComp \in \decSetFeas$ is a \textup{$\constPPE$-practically properly efficient} solution if it is properly efficient (as defined in \cref{defi:proper_efficiency}) with $\constPE \leq \constPPE$.
\end{bThm}
We can now define the term non-extreme \acc{IM}:
\begin{bThm}[Non-extreme individual minimum]{definition}
	Let $\constPPE >0$.
	We call a $\constPPE$-practically properly efficient solution $\ne{\xMO}_{\ii}^*$ a \textup{non-extreme \acc{IM}} of the $\ii$-th objective 
	if there exists no $\constPPE$-practically properly efficient solution 
	$\xMO \in \decSetFeas \setminus \ne{\xMO}_{\ii}^*$
	such that $J_{\ii}(\xMO) < J_{\ii}(\ne{\xMO}_{\ii}^*)$.
\end{bThm}
Note that we used a \form{\emph{\neAcc} accent} in order to differentiate between the 
\emph{standard} \acc{IM} $\stan{\xMO}_{\ii}^*$
and the 
\emph{non-extreme} \acc{IM} $\ne{\xMO}_{\ii}^*$.
We apply this notation to those quantities that exist in the standard and non-extreme case.
An example of this notation is the non-extreme pay-off matrix
\begin{equation}
	\ne{\Phi} := 
	\begin{mat}
		\J\(\ne{\xMO}_1^*\),&\ldots,&\J\(\ne{\xMO}_{n_J}^*\)
	\end{mat}
\end{equation}
which serves as the basis for deriving further quantities, such as the \acc{NP} and \acc{UP}.

\subsection{Distance-Based Knee-Point}
Based on the \acc{IM} and the hyperplane spanned by their convex hull, \cite{das_characterizing_1999} defines the \acc{KP} as the point in the feasible image set $\J(\decSetFeas)$ that maximizes the distance to this plane.
Furthermore, in \cite{das_normal-boundary_1998,das_characterizing_1999} it was proved that the \acc{KP} can be determined using either the \acc{WS} scalarization with suitable weights.
Since the \acc{KP} is a point that is furthest to a hyperplane, this weight vector is the normal vector $\eta$ of the considered hyperplane.
However, the normal vector is only determined up to a non-zero constant.
Since we want to avoid negative components in the weight vector, we define the following scaling operator.
\begin{bThm}[Scaling vectors]{definition}
	\label{conv:scaling_normal}
	Let the scaling operator
	$\operatorname{scal}: \Rset^n \setminus \{0\} \mapsto \Rset^n$ be defined such that for any $v \in \Rset^n \setminus \{0\}$, the output $\tilde{v} = \operatorname{scal}(v)$ is the unique vector $\tilde{v} = c v$
	with $c \in \Rset \setminus \{0\}$ chosen to satisfy
	$\operatorname{sum}(|\tilde{v}|) = 1$
	and
	$\operatorname{sign}(\tilde{v}_{i_{\max}}) = -1$, where $i_{\max} = \arg\max_{\ii} |\tilde{v}_{\ii}|$.
\end{bThm}
Note that if the normal vector $\eta$ has components with the same sign then the weight vector $w_{\mathrm{knee}}=\operatorname{scal}(\eta)$ has only non-negative components.
Similar to $\stan{w}_{\mathrm{knee}}$ based on the standard \acc{IM}, we can compute the $\ne{w}_{\mathrm{knee}}$ based on the non-extreme \accp{IM}.

\section{Computing Non-Extreme Individual Minima}
\label{sec:Derivation_Method}

\subsection{General Idea}
The development of the method starts with recapitulating one approach to obtain a \emph{standard} $\ii$-th \acc{IM}.
The approach uses the \acc{PS} scalarization \cite[Chapter 2.1]{eichfelder_adaptive_2008}
\begin{equation}
	\begin{aligned}
		\min_{\xMO \in \decSetFeas,\, l \in \mathbb{R}}
		\ & -l 
		\\
		\text{s.t.} \quad 
		& \param{\J}_{\SO} + l\param{d} - \tilde{\J}(\xMO) \in \mathsf{K}
	\end{aligned}
	\label{eq:MO_PS_standard}
\end{equation}
that is parameterized in the \acc{SO} $\param{\J}_{\SO}$ and the \acc{SDV} $\param{d}$ 
and for which we choose
\begin{subequations}
	\begin{align}
		\mathsf{K} &= 
		\left\{
			\nu \in \Rset^{n_J-1}
			\ \middle| \ 
			- \param{\spanVecMat} \nu
		\right\}
		\\
		\text{and} \quad
		\tilde{\J}(\xMO) &= \param{T}_{\J}(\J(\xMO) - \param{\J}_{\mathrm{shift}}).
	\end{align}
\end{subequations}
The parameters $\param{T}_{\J}$ and $\param{\J}_{\mathrm{shift}}$ can be used to effectively transform and shift the image space.
By appropriately setting these parameters the image space can be normalized.

In view of the different roles of the variables in this derivation, we adhere to a strict notation convention. 
While the \textit{tilde} ($\tilde{\cdot}$) is used flexibly for general modifications and transformations, the following three accents are reserved for specific contexts: 
the \textit{\paramAccent} ($\param{\cdot}$) for parameter variables, 
the \textit{\realizeAccent} ($\realize{\cdot}$) for their specific realizations, 
and the \textit{\neAcc} ($\ne{\cdot}$) for quantities associated with the non-extreme case. 
\Cref{tab:notation_accents} provides a summary of these conventions.
\begin{table}[ht]
    \centering
    \footnotesize
    \begin{tabular}{@{}lc p{8cm}@{}}
        \toprule
        Accent & Example & Reserved meaning / use case \\
        \midrule
        \textit{\paramAccent} & $\hat{p}$ & parameter of an \acc{OP} \\
        \textit{\realizeAccent} & $\realize{p}$ & numerical value of a parameter \\ 
        \textit{\neAcc}  & $\ne{\Phi}$& non-extreme case: quantities related to the non-extreme \accp[l]{IM} \\
        \textit{tilde}  & $\tilde{\mathcal{J}}$ & general modifications of a specific quantity \\
        \bottomrule
    \end{tabular}
    \vspace{0.5\baselineskip}
    \caption{Summary of mathematical notation conventions and reserved accents.}%
    \label{tab:notation_accents}%
\end{table}

We solve \eqref{eq:MO_PS_standard} with the parameter realizations
\begin{subequations}
	\label{eq:IM_param_realization}
	\begin{align}
		\realize{T}_{\J} &= I_{n_J},
		\quad
		\realize{\J}_{\mathrm{shift}} = \boldsymbol{0}
		\numberthis
		\label{eq:IM_param_realization_a}
		\\
		\realize{\J}_{\SO} &\in \Rset^{n_J},
		\quad
		\realize{d} = -e_{\ii}, 
		\quad
		\realize{\spanVecMat} = \stan{S}_{\ii} := 
			\left[e_{\kk}\right]_{\kk\in\{1,\ldots,n_J\}\setminus \ii}
		\numberthis
		\label{eq:IM_param_realization_b}
	\end{align}
\end{subequations}
to obtain the (standard) $\ii$-th \acc{IM}.
Equivalently, we can use the formulation
\begin{equation}
	\begin{aligned}
		\min_{\xMO \in \decSetFeas,\, l \in \mathbb{R}, \, \nu \in \Rset^{n_J-1} }
		\ & -l
		\\
		\text{s.t.} \quad 
		& \tilde{\J}(\xMO) = \param{\J}_{\SO} + l \param{d} + \param{\spanVecMat}\nu.
	\end{aligned}
	\tag{$P_{\mathrm{PS}}(p)$}
	\label{eq:MO_PS}
\end{equation}
with 
$
	p = 
	\{
		\param{\J}_{\SO}, 
		\param{d}, 
		\param{\spanVecMat},
		\param{T}_{\J}, 
		\param{\J}_{\mathrm{shift}},
	\}
$%
.
The \acc{OP} \eqref{eq:MO_PS} with the parameter realizations \eqref{eq:IM_param_realization} can be understood and visualized as follows:
A hyperplane, spanned by the column vectors of $\stan{S}_{\ii}$ (which we refer to as spanning vectors), is attached at the end of the shooting ray $\realize{\J}_{\SO} - l e_{\ii}$.
A feasible objective vector $\tilde{\J}(\xMO)$ must then lie on that hyperplane.
\\
We note that \eqref{eq:MO_PS} with $\realize{d}=-e_{\ii}$ and $\realize{\spanVecMat}=S_{\ii}$ 
is formally equivalent to the \acc{WS} scalarization%
\footnote{
	This formulation also covers the case of using $\tilde{\J}$ instead of $\J$:
	\\
	$
	\param{w}\tr \tilde{\J} 
	= \param{w}\tr \(\param{T}_{\J}(\J - \param{\J}_{\mathrm{shift}})\)
	= 
		\left.
			\underbrace{
				\(\param{T}_{\J}\tr \param{w}\)
			}_{\textstyle =: \param{w}^{\prime}}
		\right.\tr \J 
		- \underbrace{
			\param{w}\tr\param{\J}_{\mathrm{shift}}
		}_{\textstyle=\text{const.}}
	$
}%
\begin{equation}
	\begin{aligned}
		\min_{\xMO \in \decSetFeas}
		\ & 
		\param{w}\tr
		\J(\xMO)
	\end{aligned}
	\tag{$P_{\mathrm{WS}}(\param{w})$}
	\label{eq:MO_WS}
\end{equation}
with $\realize{w}=e_{\ii}$.
However, the application of the \acc{PS} scalarization provides a more convenient framework for the following derivation.

Instead of using basis vectors to span the hyperplane, we rotate each basis vector around a specific axis in a specific direction, see the exemplary rotation of $e_2$ in \cref{fig:general_{\ii}dea}.
\begin{figure}[tbp]
	\centering
		\adjincludegraphics[
			width=0.75\linewidth, clip, page=1,
			trim={{0\width} {0\height} {0\width} {0\height}}
		]{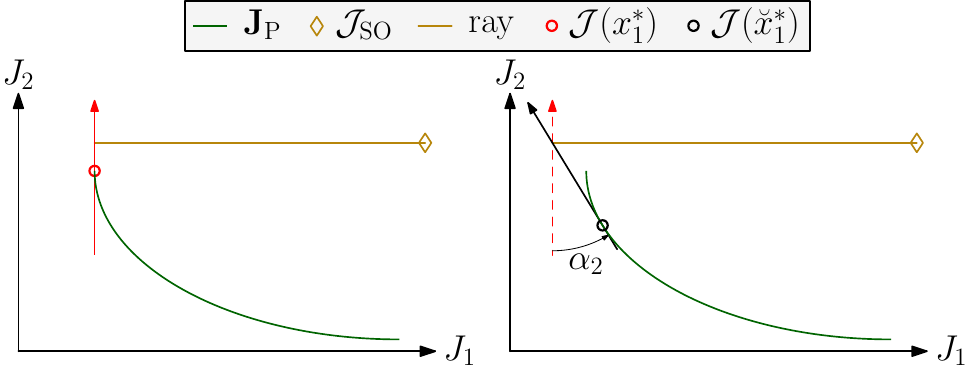}
	\caption{
		First \acc[l]{IM}: standard approach (left) and non-extreme approach (right). 
		The arrows represent the (elongated) vector $S_2$ (red) and $\ne{S}_2(\alpha)$ (black).
	}
	\label{fig:general_{\ii}dea}
\end{figure}
The rotated spanning vectors read
\begin{equation}
	\begin{aligned}
		v^{(\ii)}(\alpha_{\kk})&=
		R_{\ll,\mm}\(
			\mathrm{sign}(\kk-\ii) \, \alpha_{\kk} 
		\) e_{\kk},
		\\[0.5\baselineskip]
		\ii &\in \{1,\ldots,n_J\}, \quad 
		\kk \in \{1,\ldots,n_J\} \setminus \ii,
	\end{aligned}
	\label{eq:v_{\ii}_alpha_k}
\end{equation}
with 
$\ll = \min\{\kk,\ii\}$, 
$\mm = \max\{\kk,\ii\}$,
$\alpha = \matInline{\alpha_{1},\dots,\alpha_{n_J}}$
and $R_{\ll,\mm}(\phi)$ is a Givens rotation in the $(\ll,\mm)$-plane \cite[chapter 5.1.8]{golub_matrix_1996}.
The Givens rotation $R_{\ll,\mm}(\phi)$ is equal to the $n_J\times n_J$ identity matrix except that the entries $(\ll,\ll)$, $(\ll,\mm)$, $(\mm,\ll)$ and $(\mm,\mm)$ are overwritten with
\begin{equation}
	\begin{mat}
		\cos(\phi) & -\sin(\phi) \\
		\sin(\phi) & \cos(\phi)
	\end{mat}.
\end{equation}
The spanning vectors of our new hyperplane are concatenated horizontally to form
\begin{equation}
	\ne{S}_{\ii}(\alpha) = 
		\left[ v^{(\ii)}(\alpha_{\kk}) \right]_{\kk\in\{1,\ldots,n_J\}\setminus \ii}
	.
	\label{eq:S_{\ii}_NE}
\end{equation} 
Note that $\ne{S}_{\ii}(\alpha)=\stan{S}_{\ii}$ when $\alpha=0$. 
Furthermore, the hyperplane associated with $\ne{S}_{\ii}(\alpha)$ has a normal vector $w^{(\ii)}$ which we can scale such that $w^{(\ii)} \geq 0$ and $\sum_{\kk} w^{(\ii)}_{\kk} = 1$.

Choosing $\alpha > 0$, we can now obtain the $\ii$-th non-extreme \acc{IM} 
with \eqref{eq:MO_PS} by setting 
\begin{equation}
	\realize{\J}_{\SO} \in \Rset^{n_J},
	\quad
	\realize{d} = -e_{\ii}, 
	\quad
	\realize{\spanVecMat} = \ne{S}_{\ii}(\alpha)
	\label{eq:IM_NE_param_realization_b}
\end{equation}
and the remaining parameters as shown \eqref{eq:IM_param_realization_a}
or with \eqref{eq:MO_WS} by setting $\realize{w}=w^{(\ii)}$.
Choosing $\alpha>0$, which yields $w^{(\ii)}>0$, guarantees a finite $L$ and thus also a non-extreme \acc{IM}.

\begin{figure}[h]
	\centering
		\adjincludegraphics[
			width=0.65\linewidth, clip, 
			trim={{0.03\width} {0.0\height} {0.12\width} {0.04\height}}
		]{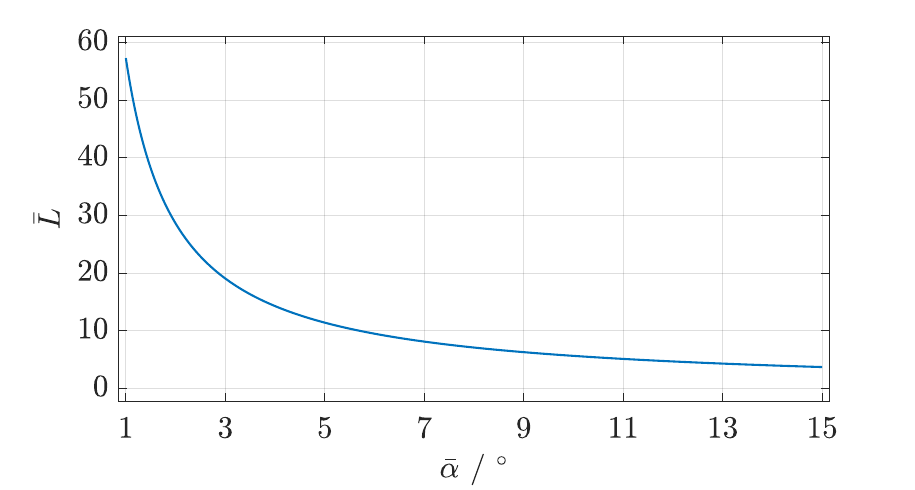}%
	\caption{Finite values of $\bar{L}$ for non-negative $\bar{\alpha}$.}
	\label{fig:L_bar_over_alpha}
\end{figure}

To see the general trend of increasing $\alpha$ components
we exemplarily set $\alpha = \bar{\alpha}$, with $\bar{\alpha} \in \Rset_{> 0}$ and compute 
\begin{equation}
 	\bar{L} := 
 	\max_{\ii \in \{1,\ldots,n_J\}} \,
 	\max_{l,m \in \{1,\ldots,n_J\}, \, l\neq m} 
 	w^{(\ii)}_l/w^{(\ii)}_m.
 	\label{eq:bar_L}
\end{equation}
As we can see in \cref{fig:L_bar_over_alpha},
with an increasing $\bar{\alpha}$ the value of $\bar{L}$ is strictly decreasing.
\form{
	Furthermore, numerical investigations show that $\bar{L}$ seems to be invariant w.r.t. the number of objectives $n_J$.
}

\subsection{Handling Different Objective Ranges}
We apply the method to the exemplary three-objective \acc{OP}
\begin{equation}
     \begin{aligned}[t]
         \min_{x \in \mathbb{R}^3} \,& \J(x) = \diag(\matInline{l_1,l_2,l_3})x
         \quad
         \text{s.t.} \; 
         x_1^2 + x_2^2 + x_3^2 \leq 1,
     \end{aligned}
     \label{eq:MOOP_ellipsoid}
\end{equation}
where the set of feasible objectives is a non-rotated ellipsoid with the semi-axis lengths $\matInline{l_1,l_2,l_3}$.
For our investigations we choose $\alpha=10^\circ$.

Our first numerical test uses the semi-axis lengths $\matInline{l_1,l_2,l_3} = \matInline{1,1,1}$ (cf.~\cref{fig:unit_sphere})
and we set the parameters as shown in \eqref{eq:IM_param_realization_a} and \eqref{eq:IM_NE_param_realization_b}.
Then, the method performs \enquote{as expected} in the sense that 
the set
\begin{equation}
	\imSetPareto[\ne] := \left\{ \J \in \imSetPareto[\stan] \mid \J \leq \ne{\J}_\NP \right\}
\end{equation}
excludes regions of the \acc{PF} that we \form{regard as practically irrelevant}.
\newlength{\SubFigWidth}%
\setlength{\SubFigWidth}{1\linewidth}%
\begin{figure}
	\centering
	\begin{subfigure}[b]{\SubFigWidth}
		\centering
			\adjincludegraphics[
				width=1\linewidth, clip, 
				trim={{0.00\width} {0.00\height} {0.04\width} {0.0\height}}
			]{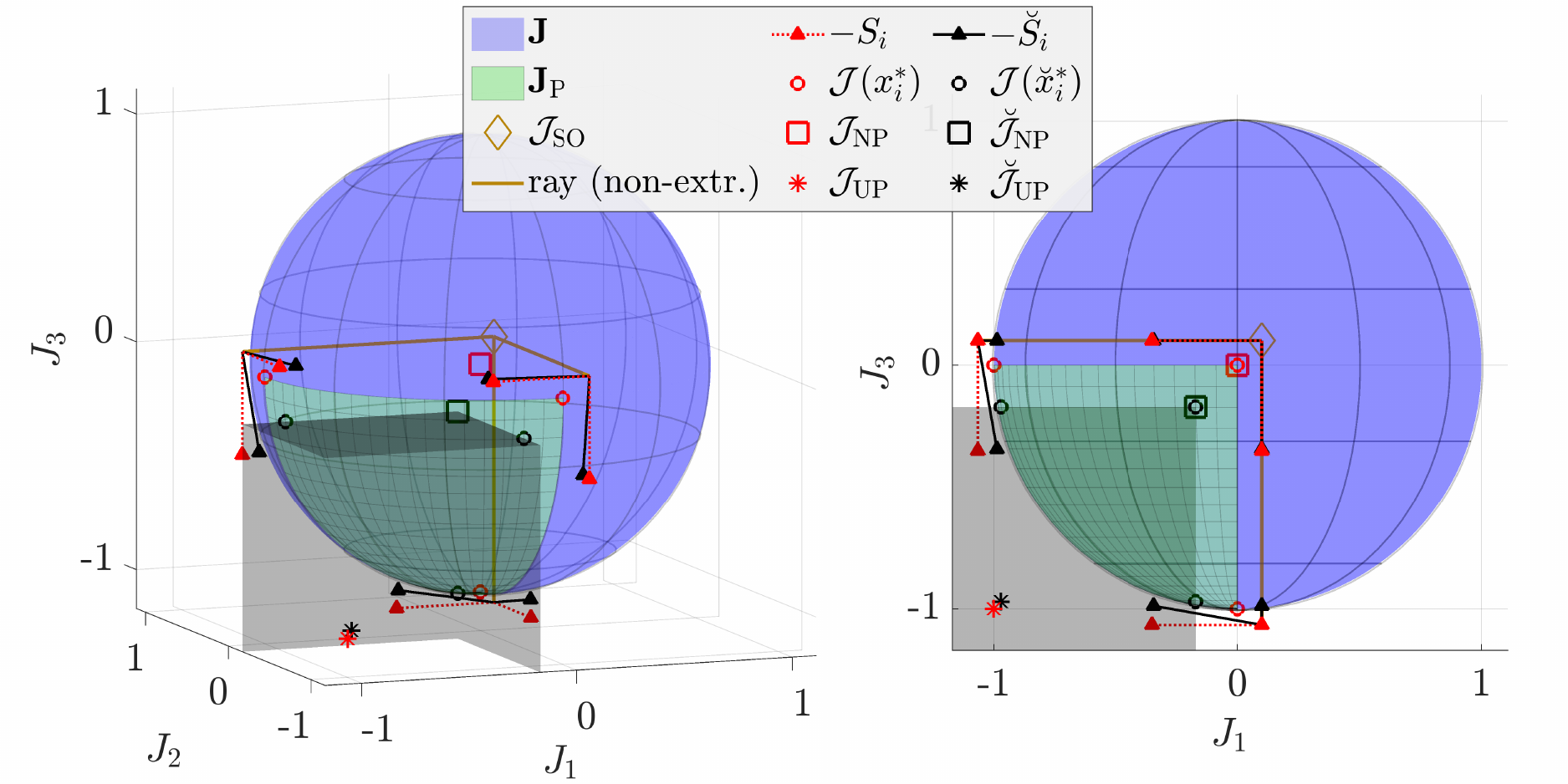}
		\caption{
			The semi-axis lengths are $\matInline{l_1,l_2,l_3}=\matInline{1,1,1}$ 
			and no image space normalization is used. 
			Note that only the ray for the non-extreme case is displayed.
			For comparison, at the end of the rays both sets of spanning vectors are \enquote{attached}.
			Due to the symmetry the views in the other two planes look the same as the $J_1$-$J_3$-plane on the right.
		}
		\label{fig:unit_sphere}
	\end{subfigure}
	\\
	\begin{subfigure}[b]{\SubFigWidth}
			\adjincludegraphics[
				width=1\linewidth, clip,
				trim={{0.00\width} {0.00\height} {0.04\width} {0.00\height}}
			]{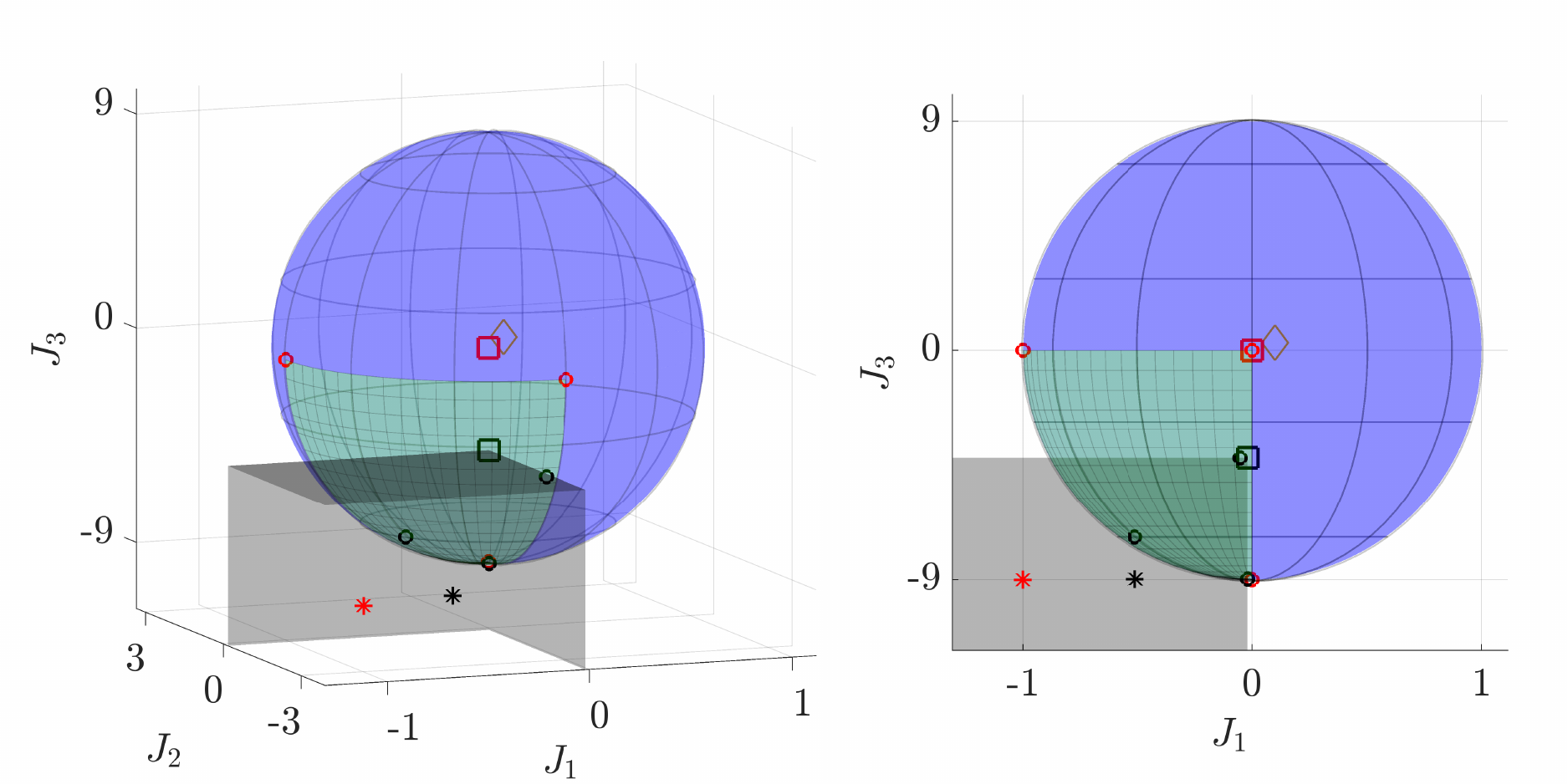}
		\caption{
			The semi-axis lengths are $\matInline{l_1,l_2,l_3} = \matInline{1,3,9}$ 
			and no image space normalization is used.
			For reasons of a cleaner visual presentation, we refrain from displaying the spanning vectors and the shooting rays.
		}
		\label{fig:ellipsoid_non_scaled}
	\end{subfigure}
	\caption{
		Results for the \acc[l]{MOOP} \protect\eqref{eq:MOOP_ellipsoid}.
	}
\end{figure}
However, in our second numerical test with $\matInline{l_1,l_2,l_3} = \matInline{1,3,9}$ (cf.~\cref{fig:ellipsoid_non_scaled}) and the parameter realizations as before
\form{the method fails to only exclude practically irrelevant regions.}

To counteract this effect we can incorporate information about the \acc{NP} and the \acc{UP} which give us valuable insight over the ranges of the objectives.
By choosing the parameter realizations
\begin{equation}
	\realize{T}_{\J} = \cNPUP,
	\quad
	\realize{\J}_{\mathrm{shift}} = \J_{\UP}
	\numberthis
	\label{eq:IM_NE_param_realization_a}
\end{equation}
we normalize the image space.
This idea of normalizing the image space also transfers to \eqref{eq:MO_WS}. 
Here, we can realize the image space normalization by setting $\realize{w}=\cNPUP\tr w^{(\ii)} = \cNPUP w^{(\ii)}$ instead of $\realize{w}=w^{(\ii)}$.

\form{%
The effect of the normalization approach is that 
(except for the shooting rays and the attached spanning vectors) 
the numerical results would produce a figure that looks like \cref{fig:unit_sphere} only that numbers on the $J_2$ and $J_3$ axis have changed to $\pm 3$ and $\pm 9$.
This means that $\imSetPareto[\ne]$ now contains the desired region of the \acc{PF}.
}

\subsection{Algorithm}
We can combine the previous findings to construct \cref{alg:NEIM}.
\begin{algorithm}[h]
	\caption{Non-Extreme Individual Minima}
	\label{alg:NEIM}
	\begin{enumerate}[label={\arabic*)},itemsep=0.7\baselineskip,labelindent=0.5ex]
		\setcounter{enumi}{-1}
		\item Use $\alpha$ to compute $\ne{S}_{\ii}(\alpha)$, 
		$\ii \in \{1,\ldots,n_J\}$. 
		Then, compute $w^{(\ii)}$ as the normal vector of the hyperplane spanned by the column vectors of $\ne{S}_{\ii}(\alpha)$
		as defined in \eqref{eq:S_{\ii}_NE} and \eqref{eq:v_{\ii}_alpha_k} 
		(and scale it as described in~\cref{conv:scaling_normal})%
		.
		\item Solve \eqref{eq:MO_WS} with $\realize{w}=e_{\ii}$ for all $\ii \in \{1,\ldots,n_J\}$ and compute $\stan{\Phi}$. 
		\item Derive $\stan{\J}_{\NP}$, $\stan{\J}_{\UP}$ and $\cNPUP$ from $\stan{\Phi}$.
		\item Solve \eqref{eq:MO_WS} with $\realize{w}=\cNPUP w^{(\ii)}$ for all $\ii \in \{1,\ldots,n_J\}$ and compute $\ne{\Phi}$.
	\end{enumerate}	
\end{algorithm}
Notably, step 0) of this algorithm needs to be performed only once for a fixed $\alpha$.
Furthermore, if the values of $\J_{\NP}$ and $\J_{\UP}$ are not available, which is generally the case, step 1) and 2) have to be executed.
Then, in total, $2 n_J$ \accp{OP} have to be solved to determine the non-extreme \accp{IM}.

\begin{figure}[h]
	\centering
		\adjincludegraphics[
			width=0.73\linewidth, clip, 
			trim={{0\width} {0\height} {0\width} {0\height}}
		]{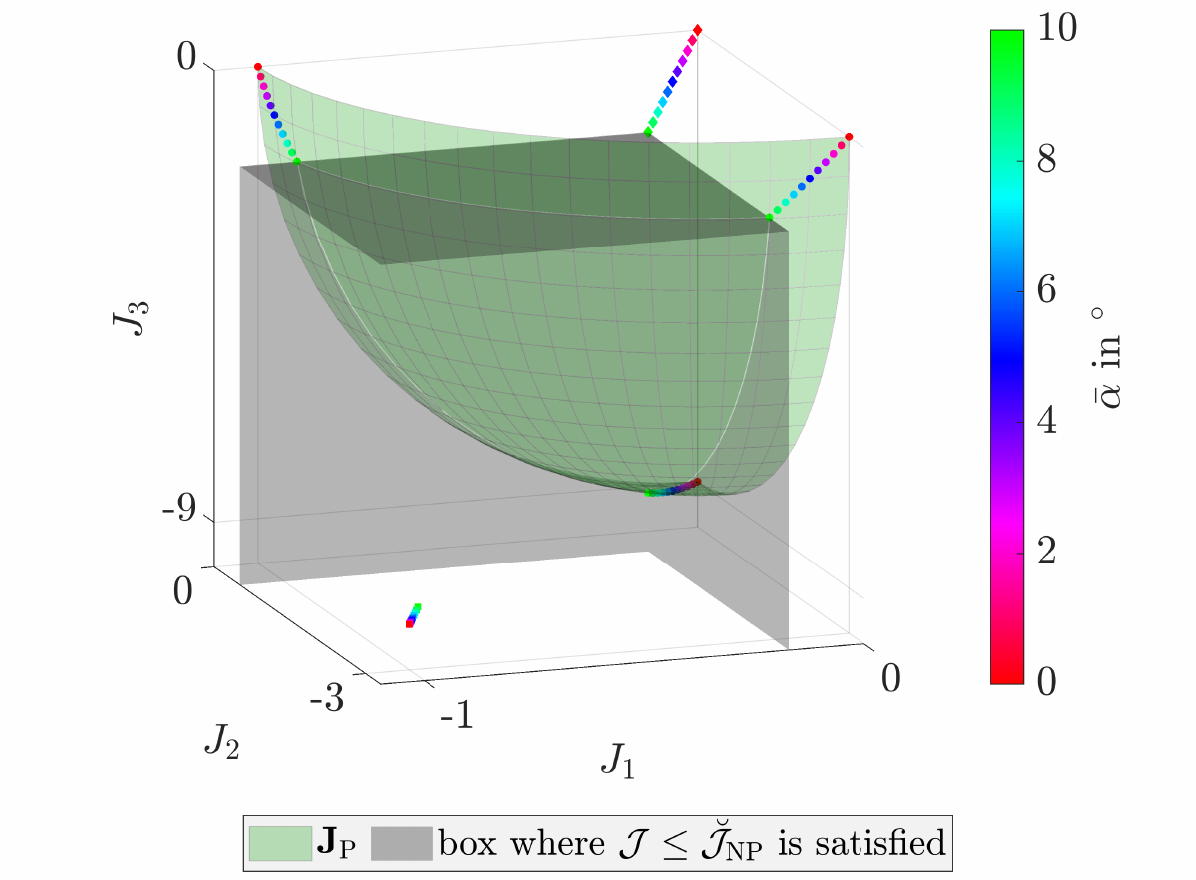}
	\caption{The colored markers are related to the non-extreme \accp{IM} (with $\alpha=\bar{\alpha}$): dots represent the non-extreme \accp{IM}, while diamonds and squares represent the corresponding \accp{NP} and \accp{UP}, respectively.}
	\label{fig:NEIM_over_alpha}
\end{figure}

Applying \cref{alg:NEIM} to \acc{MOOP} \eqref{eq:MOOP_ellipsoid} with $\alpha=\bar{\alpha}$, where $\bar{\alpha}=0^\circ,1^\circ,\ldots,10^\circ$, yields \cref{fig:NEIM_over_alpha}.
While \acc{NP} values of non-extreme \accp{IM} drop significantly, the associated \acc{UP} values show only a slight increase.
This effect yields hypercubes of reduced size which allow for 
	a denser sampling of the \acc{PF} with a fixed number of sampling points 
	or, alternatively, fewer samples for a targeted approximate point density.
By omitting practically irrelevant regions, improved resolution or efficiency gains are achieved.

\section{Numerical Example}
\label{sec:Numerical_Example}

\begin{figure*}[tbp]
	\centering
		\adjincludegraphics[
			width=1\linewidth, clip,
			trim={{0.025\width} {0.00\height} {0.06\width} {0.00\height}}
		]{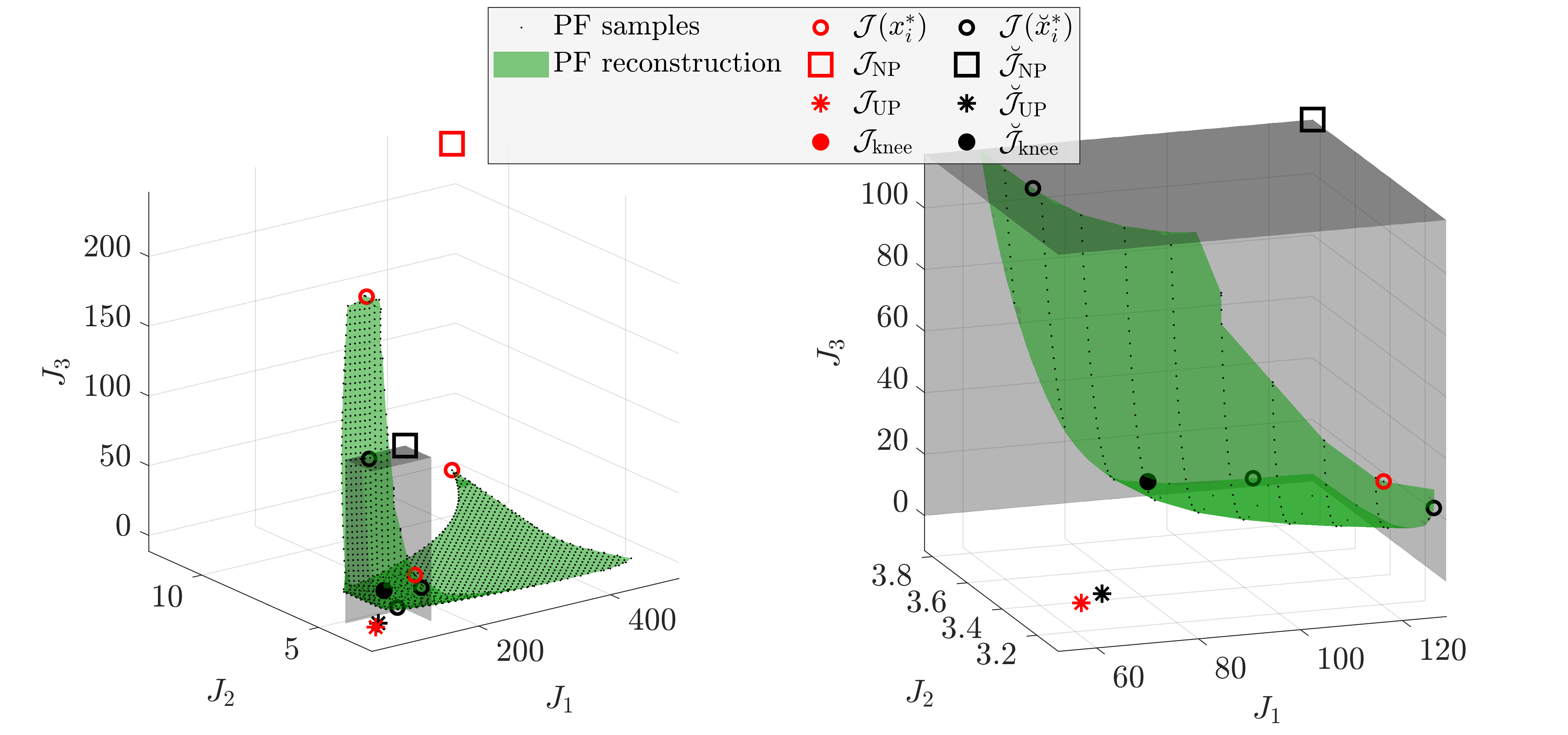}%
	\caption{
		\acc{PF} reconstruction using the \acc{PF} samples and the ball-pivoting algorithm~\cite{bernardini_ball-pivoting_1999}. 
		The right plot shows a segment of the \acc{PF} that does not contain \form{practically irrelevant} regions.
		Both figures use an equal plot box aspect ratio which has a similar visual effect as normalizing the image space and allows for a fair comparison of both \acc{PF} regions.
	}
	\label{fig:HVAC_PF_NEIM}
\end{figure*}

As our numerical example we consider the \acc{PF} of \enquote{Figure 6.} from \cite{herrmann-wicklmayr_individual_2025}.
This reference describes in detail how the \acc{MOOP}, 
which 
	implements a nonlinear \emph{\acc{HVAC}} control problem
	and yields the \acc{PF} depicted in \cref{fig:HVAC_PF_NEIM} (left), 
is constructed.
We choose this specific \acc{PF} because i) it constitutes a non-academic example and ii) it has a large region where it is either extremely flat or steep.

We set $\alpha=3^{\circ}$ which results in $L=\bar{L}<20$ (cf. \eqref{eq:bar_L}) for the non-extreme \accp{IM} (cf.~\cref{fig:L_bar_over_alpha}).
We then compute the standard \accp{IM} and the non-extreme \accp{IM} as described in \cref{alg:NEIM}.
Both resulting sets of points can be seen in \cref{fig:HVAC_PF_NEIM}.
As can be concluded from the closeness of $\ne{\J}_{\UP}$ to $\stan{\J}_{\UP}$ only extreme trade-off solutions are excluded.
This is supported by the fact that \cref{fig:HVAC_PF_NEIM} (right) has no regions that are extremely flat or steep.

We note that from all 1520 \acc{PF} samples shown in \cref{fig:HVAC_PF_NEIM} (left) only 186 are contained in the box spanned by $\ne{\J}_\NP$ and $\ne{\J}_\UP$.
This means that, since sampling \form{practically irrelevant} points on the \acc{PF} provides no added value to the decision-maker, more than 87\% of the optimizations represented a superfluous expenditure of time and resources that our method successfully eliminates.
Furthermore, for \cref{fig:HVAC_PF_NEIM}, the \accp{KP} based on the standard and non-extreme \accp{IM} were computed. 
This requires the normal vector of the convex hull of the \accp{IM} and non-extreme \accp{IM}, which are $\stan{w}_{\mathrm{knee}} = \matInline{-0.0239, 0.9673, -0.0088}\tr$ and $\ne{w}_{\mathrm{knee}} = \matInline{0.0306, 0.9527, 0.0167}\tr$, respectively.
Note that $\stan{w}_{\mathrm{knee}}$ contains negative components. 
This leads to two critical issues: 
first, the theoretical guarantee of obtaining a Pareto optimal solution is lost; 
second, and more importantly, the scalarization effectively rewards increases in cost for objectives associated with negative weights, contradicting the fundamental goal of minimization.
Consequently, the resulting point $\stan{\J}_{\mathrm{knee}} = 10^3 \cdot \matInline{2.9783, 0.0036, 3.1886}\tr$ is a dominated point located on the boundary of the feasible image set $\J(\decSetFeas)$ and is outside the axis intervals of \cref{fig:HVAC_PF_NEIM}.
Enforcing a lower bound of zero on the weight components as a \enquote{safety layer} would, in this example, lead to the recovery of the second \acc{IM}. 
However, such a solution contradicts the fundamental concept of a knee point, which is intended to represent a balanced compromise rather than an extreme boundary solution.
In contrast, the solution $\ne{\J}_{\mathrm{knee}}$ obtained by using the weight $\ne{w}_{\mathrm{knee}}$ 
(constructed from the non-extreme \accp{IM}) 
constitutes a balanced trade-off.

\section{Conclusion and Outlook}
\label{sec:Conclusion_and_Outlook}
In this paper, we introduced the concept of non-extreme \accp[l]{IM} as a means to exclude \form{practically irrelevant} regions of the \acc[l]{PF}.
By leveraging the definition of $L$-practical proper efficiency, we derived a method to exclude solution candidates that represent unreasonable trade-offs.

The proposed algorithm is straightforward to implement, relying solely on \acc[l]{WS} scalarizations with modified weight vectors obtained via geometric rotations.
We demonstrated that integrating information about the \acc[l]{UP} and the \acc[l]{NP} to normalize the image space is crucial for the method's robustness against differing objective ranges.
The numerical results confirm that the computed non-extreme \acc[l]{IM} effectively bound the area of interest, excluding extremely steep or flat parts of the front.
This bounding box allows subsequent multi-objective optimization methods or automated decision-making schemes to focus their computational effort on the most promising trade-offs.

Future work could evaluate the influence of this approach on algorithms that rely on automated \acc[l]{DM}. 
A prime example of this is the \acc[l]{PF} from \cite{herrmann-wicklmayr_individual_2025} shown in \cref{fig:HVAC_PF_NEIM}, which arises during the first iteration of an \acc{MOMPC} scheme.

\bibliography{25_Non_Extreme_IM.bib}

\end{document}